\documentclass{article}
\usepackage{graphicx} 
\usepackage{placeins}
\usepackage{amssymb}
\usepackage{amsmath}
\usepackage{amsfonts}
\usepackage{amsthm}
\usepackage{enumitem}
\usepackage{hyperref}
\usepackage[dvipsnames]{xcolor}

\usepackage{tikz}
\usetikzlibrary{matrix}
\usetikzlibrary{positioning}
\usetikzlibrary {patterns,patterns.meta}
\usetikzlibrary{shapes.geometric}

\theoremstyle{remark}
\newtheorem{remark}{Remark}

\theoremstyle{definition}
\newtheorem*{ex}{Exercise}
\newtheorem{defi}{Definition}
\newtheorem{theorem}{Theorem}
\newtheorem{lemma}{Lemma}

\newcommand{\R}{\mathbb{R}}
\newcommand{\N}{\mathbb{N}}

\title{Constructing many-twist M\"obius bands with small aspect ratios}
\author{Aidan Hennessey}
\date{April 2024}

\begin{document}

\maketitle

\begin{abstract}
    This paper presents a construction of a folded paper ribbon knot that provides a constant upper bound on the infimal aspect ratio for paper M\"obius bands and annuli with arbitrarily many half-twists. In particular, the construction shows that paper M\"obius bands and annuli with any number of half-twists can be embedded with aspect ratio less than 8.
\end{abstract}

\section{Introduction}

In 1977, Halpern and Weaver conjectured that the infimal aspect ratio of an embedded paper M\"obius band is $\sqrt{3}$ \cite{hw}. This conjecture was recently proven by Schwartz in \cite{one-twist}. Shortly after, Schwartz proved that the infimal aspect ratio for the 2 half-twist annulus is 2 \cite{two-twist}. Noah Montgomery independently showed this result using alternative methods (unpublished). Moreover, Brown and Schwartz conjecture the infimal aspect ratio for the 3-half-twist embedded paper M\"obius band is 3 \cite{three-twist}. Beyond these published results, unpublished experiments by Brown indicate the infimal aspect ratios for embedded paper M\"obius bands with 5 and 7 half-twists are each at most 5.

The seeming pattern of more twists requiring a longer band raises the question: What is the asymptotic growth of the minimal aspect ratio $\lambda_n$ as a function of the number of twists $n$? Noah Montgomery found a construction with length complexity $O(\sqrt{n})$ (unpublished), but did not produce any lower bound. This paper puts the big-$O$ question to rest by constructing an $O(1)$ solution. The paper does not show that the construction's constant bound is tight. I.e., determining the value of $\lim \sup \{\lambda_n\}$ is still an open problem. The construction is actually a folded ribbon (un)knot which can be arbitrarily well-approximated by paper bands. Its folded ribbon knot form negatively answers Conjecture 39 in \cite{ribbon-knots}.

\begin{theorem}[Main Theorem]
    There exists a constant $\lambda$ such that for any $n$, there exists a paper band (defined below) with $n$ half-twists and aspect ratio less than $\lambda$. In fact, it suffices to let $\lambda = 8$.
\end{theorem}

This theorem is an immediate corollary of two lemmas. The first lemma is a statement about a family of objects known as folded ribbon knots. Roughly speaking, a folded ribbon knot is a folded strip of paper which lies in the plane (See Remark 1 for more). The \emph{ribbon linking number} of a folded ribbon knot $\mathcal{K}$ is linking number between its centerline and one boundary component. The folded ribbon length $\text{Rib}(\mathcal{K})$ of $\mathcal{K}$ is the aspect ratio of the strip of paper, before folding. See Definition 3 for more detail.

\begin{lemma}[High-Link Paper Ribbon Knots]
    There is a family $\{\mathcal{K}_n\}$ of folded ribbon knots such that 
    \begin{itemize}[topsep=0pt, itemsep=0pt]
        \item If $n$ is odd, $\mathcal{K}_n$ is a topological M\"obius band with ribbon linking number $\pm n$.
        \item If $n$ is even, $\mathcal{K}_n$ is a topological annulus with ribbon linking number $\pm n/2$.
        \item There exists a constant $\lambda$ such that for all $n \in \N$, $\mathcal{K}_n$ has folded ribbon length $\text{Rib}(\mathcal{K}_n) < \lambda$. In fact, $\lambda = 8$ suffices.
    \end{itemize}
\end{lemma}

\begin{lemma}[Approximability]
    For each $n$, there is a sequence of $n$-twist paper bands in $\R^3$ which converge pointwise to the folded ribbon knot $\mathcal{K}_n$, and whose aspect ratios converge to $\text{Rib}(\mathcal{K}_n)$. 
\end{lemma}

Section 2 introduces necessary definitions and terminology. 
Section 3.1 defines a particular family $\{\mathcal{K}_n\}_{n \in \N}$, and shows that it satisfies bullet 3 of Lemma 1. Section 3.2 shows that this family satisfies bullets 1 and 2, completing the proof of Lemma 1. Section 3.3 proves Lemma 2. Section 4 gives an explicit value for the bounding aspect ratio $\lambda$ of Theorem 1.

\section{Background}

\begin{defi}[Paper Band, Aspect Ratio]
    Formally, a paper M\"obius band is a smooth locally isometric embedding of the M\"obius strip $$([0, 1] \times [0, \lambda]) / \sim; \>\>\>\>\>\>\>(t, 0) \sim (1-t, \lambda)$$ into $\mathbb{R}^3$. Similarly, a paper annulus is a smooth locally isometric embedding of the cylinder $([0, 1] \times [0, \lambda]) / ((t, 0) \sim (t, \lambda))$ into $\mathbb{R}^3$. Refer to these maps collectively as \emph{paper bands}. $\lambda$ is called the \emph{aspect ratio} of the band.
\end{defi}

\begin{defi}[Center line, half twist]
    The \textit{center line} of a band is the image of $\{0.5\} \times [0, \lambda]$ under the embedding. Define an $n$ \emph{half-twist} paper band to be
    \begin{itemize}[itemsep=0pt, topsep=6pt]
        \item A paper M\"obius band for which the boundary and center line have linking number $\pm n$ (for $n$ odd).
        \item A paper annulus for which one of the boundaries and the center line have linking number $\pm n/2$ (for $n$ even).
    \end{itemize}
\end{defi}

\begin{remark}
    Many\footnote{It is likely not all paper bands have this property. A particular likely counterexample is the cap, an efficient 3-twist band featured in \cite{three-twist}. This counterexample was pointed out to me by Richard Schwartz.} paper bands can be gently pressed down to lie in the plane, at which point the image is a union of rectangles, parallelograms, and trapezoids, joined to one another at creases in sequence. Such an object is known as a \emph{folded ribbon knot}. For a formal definition, see \cite{ribbon-knots}. 
\end{remark}

\begin{defi}[Folded Ribbon Length]
    The centerline of a folded ribbon knot $\mathcal{K}$ is a closed polygonal curve. The ratio of the total length of this curve to the ribbon knot's width is the \emph{folded ribbon length} of $\mathcal{K}$, denoted $\text{Rib}(\mathcal{K})$. In this paper, all widths are taken to be 1, so this is just the length of the centerline.
\end{defi}

\begin{defi}[Prefold Diagram]
    To a folded ribbon knot we can associate a \emph{prefold diagram}, which is a rectangle with non-intersecting solid and dotted line segments (prefolds) drawn on it. Each line segment represents a fold, and the texture of the segment dictates which way the fold goes. One can imagine the rectangle as a strip of paper, with the side facing the viewer colored red, and the other side colored blue. Then, a solid line indicates folding so that the red side is on the inside, and a dotted line indicates a fold which has a blue inside.
\end{defi}

\vspace{-5pt}

\begin{figure}[!hb]
    \centering
    \includegraphics[scale=0.3]{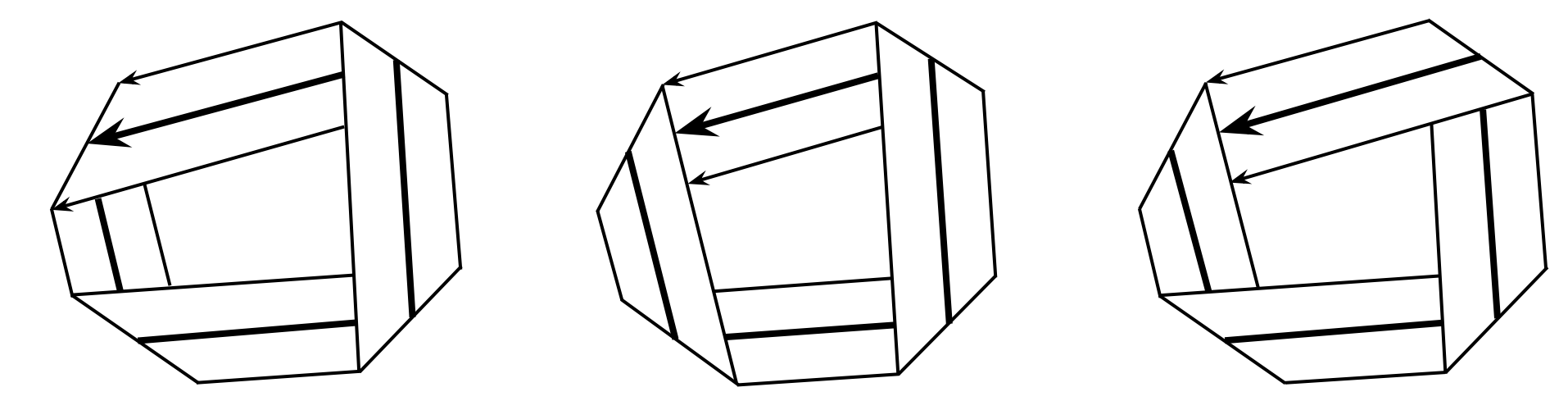}
    \caption{Three different ribbon knots. Their prefold diagrams all have line segments in the same places, but they differ in which segments are solid or dashed. Reused with permission from \cite{figures}.}
    \label{fig:borrowed}
\end{figure}

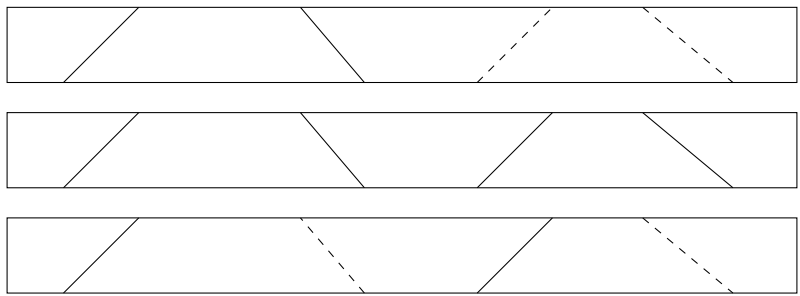
\begin{figure}[!h]

    \centering
    
    \begin{tikzpicture}
        \node[rectangle, minimum height=1cm, minimum width=10.5cm, draw] {};
        \draw (-4.5, -0.5)--(-3.5, 0.5);
        \draw (-0.5, -0.5)--(-1.35, 0.5);
        \draw[dashed] (1, -0.5)--(2, 0.5);
        \draw[dashed] (3.2, 0.5)--(4.4, -0.5);
    \end{tikzpicture}

    \vspace{10pt}
    
    \begin{tikzpicture}
        \node[rectangle, minimum height=1cm, minimum width=10.5cm, draw] {};
        \draw (-4.5, -0.5)--(-3.5, 0.5);
        \draw (-0.5, -0.5)--(-1.35, 0.5);
        \draw (1, -0.5)--(2, 0.5);
        \draw (3.2, 0.5)--(4.4, -0.5);
    \end{tikzpicture}
    
    \vspace{10pt}
    
    \begin{tikzpicture}
        \node[rectangle, minimum height=1cm, minimum width=10.5cm, draw] {};
        \draw (-4.5, -0.5)--(-3.5, 0.5);
        \draw[dashed] (-0.5, -0.5)--(-1.35, 0.5);
        \draw (1, -0.5)--(2, 0.5);
        \draw[dashed] (3.2, 0.5)--(4.4, -0.5);
    \end{tikzpicture}
    \caption{The prefold diagrams for the above three ribbon knots. The top prefold diagram corresponds to the left ribbon knot, the middle with middle, and bottom with right.}
    \label{fig:prefold-ex}
\end{figure}


\section{Construction}

\subsection{Folded Ribbon Knot}

Most constructions aimed at this problem are centered around the ``belt trick": Coil a belt, and then pull the ends apart without allowing them to rotate. The coils turn into twists. This is useful because it means one can construct a many-twist band by tightly coiling the band, yielding very many twists while using a small length of band. The issue with this is that one end of the band ends up confined in a very small space, which prevents reconnection of the two ends without using a very large amount of band to ``escape." 

Here's the key idea for this paper: If we wrap very tightly at a large angle, then we can escape using a constant length of band. Construct an ``escape accordion" by folding along parallel lines, 45 degrees rotated from the sides of the band. 

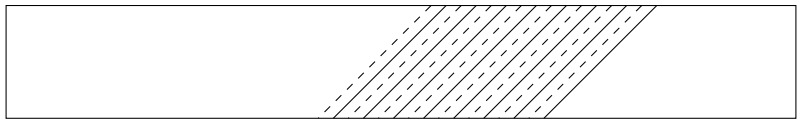
\begin{figure}[!h]
    \centering
    \begin{tikzpicture}

    \node[rectangle, draw, minimum height=1.5cm, minimum width=10.5cm] at (-2, 0) {}; 

    \draw (1.4, 0.75)--(-0.1, -0.75);
    \draw[dashed] (1.2, 0.75)--(-0.3, -0.75);
    \draw (1.0, 0.75)--(-0.5, -0.75);
    \draw[dashed] (0.8, 0.75)--(-0.7, -0.75);
    \draw (0.6, 0.75)--(-0.9, -0.75);
    \draw[dashed] (0.4, 0.75)--(-1.1, -0.75);
    \draw (0.2, 0.75)--(-1.3, -0.75);
    \draw[dashed] (0, 0.75)--(-1.5, -0.75);
    \draw (-0.2, 0.75)--(-1.7, -0.75);
    \draw[dashed] (-0.4, 0.75)--(-1.9, -0.75);
    \draw (-0.6, 0.75)--(-2.1, -0.75);
    \draw[dashed] (-0.8, 0.75)--(-2.3, -0.75);
    \draw (-1, 0.75)--(-2.5, -0.75);
    \draw[dashed] (-1.2, 0.75)--(-2.7, -0.75);
    \draw (-1.4, 0.75)--(-2.9, -0.75);
    \draw [dashed] (-1.6, 0.75)--(-3.1, -0.75);
        
    \end{tikzpicture}
    \caption{The prefold diagram for the escape accordion}
    \label{fig:accordion-pre}
\end{figure}

Color the front side of the band red and the back side blue. Then, after folding the accordion, the band looks like this:

\begin{figure}[!h]
    \vspace{60pt}
    \hspace{110pt}
    \begin{tikzpicture}[transform canvas={scale=0.5}]
        \node[trapezium, trapezium left angle=90, trapezium right angle=135, draw, fill=red!10, minimum height=2cm, minimum width=6cm] (bottom) at (0, 0) {};

        \node[trapezium, trapezium left angle=45, trapezium right angle=90, draw, fill=red!10, minimum height=2cm, minimum width=6cm] (top) at (8.5, 3) {};

        \node[isosceles triangle, isosceles triangle apex angle=90, rotate=-45, draw, fill=blue!10] at (4.33, 1.66) {};
        \node[isosceles triangle, isosceles triangle apex angle=90, rotate=-45, draw, fill=blue!10] at (4.83, 2.16) {};
        \node[isosceles triangle, isosceles triangle apex angle=90, rotate=-45, draw, fill=blue!10] at (5.33, 2.66) {};
        \node[isosceles triangle, isosceles triangle apex angle=90, rotate=-45, draw, fill=blue!10] at (5.83, 3.16) {};
        \node[isosceles triangle, isosceles triangle apex angle=90, rotate=-45, draw, fill=blue!10] at (6.33, 3.66) {};
        \node[isosceles triangle, isosceles triangle apex angle=90, rotate=-45, draw, fill=blue!10] at (6.83, 4.16) {};

        \node[isosceles triangle, isosceles triangle apex angle=90, rotate=135, draw, fill=red!10] at (4.67, 1.82) {};
        \node[isosceles triangle, isosceles triangle apex angle=90, rotate=135, draw, fill=red!10] at (5.17, 2.32) {};
        \node[isosceles triangle, isosceles triangle apex angle=90, rotate=135, draw, fill=red!10] at (5.67, 2.82) {};
        \node[isosceles triangle, isosceles triangle apex angle=90, rotate=135, draw, fill=red!10] at (6.17, 3.32) {};
        \node[isosceles triangle, isosceles triangle apex angle=90, rotate=135, draw, fill=red!10] at (6.67, 3.82) {};
    \end{tikzpicture}
    \vspace{20pt}
    \caption{The escape accordion made from colored paper}
    \label{fig:accordion-post}
\end{figure}
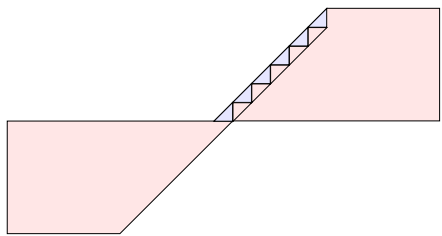

The key insight about the accordion is that its construction uses a parallelogram with base 2,\footnote{Recall that here and throughout the paper, bands and ribbon knots are assumed to have width 1.} regardless of the distance between folds. Thus, if we want to achieve $n$ half-twists using $\epsilon$ additional length of band, we can let there be $\lceil n/2\epsilon \rceil$ folds in the accordion. A base of 2 is required so that the two ends of the band do not crash into each other during the wrapping step (see Fig. 6). Adding in the prefolds for the wrapping step yields a new prefold diagram:

\newpage

\begin{figure}[ht!]
    \vspace{-15pt}
    \centering
    \begin{tikzpicture}

    \node[rectangle, draw, minimum height=1.5cm, minimum width=10.5cm] at (-2, 0) {}; 

    \draw (1.4, 0.75)--(-0.1, -0.75);
    \draw[dashed] (1.2, 0.75)--(-0.3, -0.75);
    \draw (1.0, 0.75)--(-0.5, -0.75);
    \draw[dashed] (0.8, 0.75)--(-0.7, -0.75);
    \draw (0.6, 0.75)--(-0.9, -0.75);
    \draw[dashed] (0.4, 0.75)--(-1.1, -0.75);
    \draw (0.2, 0.75)--(-1.3, -0.75);
    \draw[dashed] (0, 0.75)--(-1.5, -0.75);
    \draw (-0.2, 0.75)--(-1.7, -0.75);
    \draw[dashed] (-0.4, 0.75)--(-1.9, -0.75);
    \draw (-0.6, 0.75)--(-2.1, -0.75);
    \draw[dashed] (-0.8, 0.75)--(-2.3, -0.75);
    \draw (-1, 0.75)--(-2.5, -0.75);
    \draw[dashed] (-1.2, 0.75)--(-2.7, -0.75);
    \draw (-1.4, 0.75)--(-2.9, -0.75);
    \draw [dashed] (-1.6, 0.75)--(-3.1, -0.75);
    \draw (-1.8, 0.75)--(-3.3, -0.75);
    \draw (-2.0, 0.75)--(-3.5, -0.75);
    \draw (-2.2, 0.75)--(-3.7, -0.75);
    \draw (-2.4, 0.75)--(-3.9, -0.75);
    \draw (-2.6, 0.75)--(-4.1, -0.75);
    \draw (-2.8, 0.75)--(-4.3, -0.75);
    \draw (-3.0, 0.75)--(-4.5, -0.75);
        
    \end{tikzpicture}
    \caption{The prefold diagram for the accordion and the wrapping}
    \label{fig:enter-label}
\end{figure}
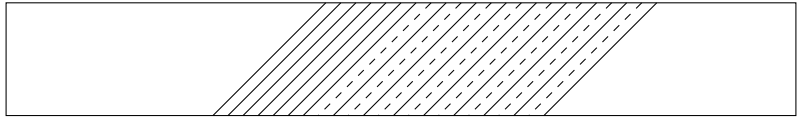

We can now fold this up, retaining the same red-blue coloring used in Figure 4, to obtain Figure 6.

\begin{figure}[!h]
    \vspace{50pt}
    \hspace{125pt}
    \begin{tikzpicture}[transform canvas={scale=0.45}]
        \node[trapezium, trapezium left angle=90, trapezium right angle=135, draw, fill=red!10, minimum height=2cm, minimum width=6cm] (bottom) at (-1.43, -1.44) {};

        \node[trapezium, trapezium left angle=45, trapezium right angle=90, draw, fill=red!10, minimum height=2cm, minimum width=6cm] (top) at (8.5, 3) {};

        \node[isosceles triangle, isosceles triangle apex angle=90, rotate=-45, draw, fill=blue!10] at (4.33, 1.66) {};
        \node[isosceles triangle, isosceles triangle apex angle=90, rotate=-45, draw, fill=blue!10] at (4.83, 2.16) {};
        \node[isosceles triangle, isosceles triangle apex angle=90, rotate=-45, draw, fill=blue!10] at (5.33, 2.66) {};
        \node[isosceles triangle, isosceles triangle apex angle=90, rotate=-45, draw, fill=blue!10] at (5.83, 3.16) {};
        \node[isosceles triangle, isosceles triangle apex angle=90, rotate=-45, draw, fill=blue!10] at (6.33, 3.66) {};
        \node[isosceles triangle, isosceles triangle apex angle=90, rotate=-45, draw, fill=blue!10] at (6.83, 4.16) {};

        \node[isosceles triangle, isosceles triangle apex angle=90, rotate=135, draw, fill=red!10] at (4.67, 1.82) {};
        \node[isosceles triangle, isosceles triangle apex angle=90, rotate=135, draw, fill=red!10] at (5.17, 2.32) {};
        \node[isosceles triangle, isosceles triangle apex angle=90, rotate=135, draw, fill=red!10] at (5.67, 2.82) {};
        \node[isosceles triangle, isosceles triangle apex angle=90, rotate=135, draw, fill=red!10] at (6.17, 3.32) {};
        \node[isosceles triangle, isosceles triangle apex angle=90, rotate=135, draw, fill=red!10] at (6.67, 3.82) {};

        \node[trapezium, fill=red!10, trapezium left angle=135, trapezium right angle=45, rotate=45, draw, minimum width=1.01cm] at (4.01, 1.26) {};
        \node[trapezium, fill=red!10, trapezium left angle=135, trapezium right angle=45, rotate=45, draw, minimum width=1.02cm] at (3.53, 0.78) {};
        \node[trapezium, fill=red!10, trapezium left angle=135, trapezium right angle=45, rotate=45, draw, minimum width=1.02cm] at (3.05, 0.30) {};

        \draw (3.1, -0.1)--(3.2, -0.2)--(4.7, 1.3)--(4.6, 1.4);
        \draw (4, 0.6)--(4.2, 0.4);
        \node at (4.6, 0) { \Huge $\delta_0$};
    \end{tikzpicture}
    \vspace{30pt}
    \caption{The full construction, up to reattaching the ends. Notice that for any fixed number of twists, the distance labeled $\delta_0$ can be made arbitrarily small with an adequately skinny accordion.}
    \label{fig:main construction}
\end{figure}
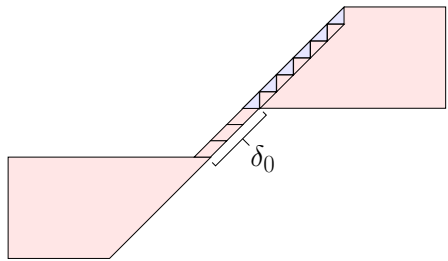

All that's left in the construction is to reattach the ends, which requires a length of band $l$ independent of the number of twists. Fixing an $\epsilon>0$ of our choosing, we recover a family $\{\mathcal{K}_n\}_{n \in \N}$ of folded ribbon knots with $\text{Rib}(\mathcal{K}_n) < 2 + l + \epsilon$.

\FloatBarrier

\subsection{Linking Number}

To prove the remainder of Lemma 1, make use of the following result:

\begin{lemma}[\cite{ribbon-knots}, Lemma 11]
    The ribbon linking number of a folded ribbon knot is determined by the combinatorial information of its folds and crossings. Each fold and crossing has a certain ``local contribution," and the sum of these local contributions is the ribbon linking number.
\end{lemma}

The term ``crossing" here refers to a place where the centerline crosses itself. General overlap of the folded ribbon knot with itself does not count as a crossing.

Here's the kicker: the knots $\mathcal{K}_n$ do not have any crossings, so the only information relevant to calculating the linking number is the folding information. There are four types of folds to consider, each with their own local contribution and realization in the prefold diagram:

\begin{enumerate}[itemsep=0pt, topsep=0pt]
    \item Right underfolds contribute $+1$ to the ribbon linking number and appear as downward sloping dashed lines in prefold diagrams
    \item Right overfolds, $-1$ - downward sloping solid lines
    \item Left underfolds, $-1$ - upward sloping dashed lines
    \item Left overfolds, $+1$ - upward sloping solid lines
\end{enumerate}

Note that the above contributions are \emph{only for M\"obius strips}. In a M\"obius strip, either side of the center line is part of the same single boundary. Compared to an annulus, then, each fold creates twice as many intersections between the centerline and a boundary component, and thus contributes twice as much to the linking number in the M\"obius band case versus the annulus case. Hence, the contribution of a fold in the annulus case is $\pm \frac{1}{2}$, not $\pm 1$. The distinct cases in the definition of an $n$ half-twist paper band exist to counterbalance this artifact.

\begin{ex}
    Using the described method, calculate the ribbon linking number of the ribbons corresponding to the prefold diagrams in Figure 2. Confirm that your answer matches what you would visually infer from Figure 1.
\end{ex}

Using this counting method, we can see in Figure 5 that the folds of the accordion cancel out in pairs, while the folds of the wrapping step compound, causing the linking number to accumulate. It takes $n$ consecutive solid lines in the prefold diagram to create a band with $n$ half-twists. Note that Figure 5 does not include the prefolds corresponding to how the ends are reconnected. The folds which are added to ensure the ends of the band connect will contribute some additional linking or unlinking, but any reasonable method only contributes a constant amount, so this does not matter. Lemma 1 is thus proven.

\subsection{Smooth Approximation}


We now prove Lemma 2, which states that the folded ribbon knots $\mathcal{K}_n$ can be well-approximated by smooth M\"obius bands.

\begin{proof}

The main idea is to model each fold with a very tight turnaround, or pseudofold. For other examples of a similar procedure, see \cite{hw, two-twist, three-twist}.

\begin{defi}[Pseudofold]
    As defined in \cite{hw}, proof of Lemma 9.1, a pseudofold is based on a plane curve $\gamma(\delta, t)$ (parameterized by arc length $t$) with curvature $\kappa(t)$ satisfying:
    \begin{itemize}[topsep=0pt, itemsep=0pt]
        \item $\kappa(t)$ is smooth and has compact support (bump function)
        \item $\kappa(t) \geq 0$
        \item $\int \kappa(t) = \pi$
    \end{itemize}

    $\gamma(\delta, t)$ follows the $x$-axis for some time, turns around smoothly, and then follows the line $y=\delta$ is the other direction. The length of the curved part is $c \delta$ for some constant $c$ depending on the particular bump function chosen. Let the curved part correspond to $t \in [0, c \delta]$

    Given such a curve $\gamma$, one can construct the chart $$(t, s) \mapsto (dt+s, \gamma_x(\delta, t), \gamma_y(\delta, t))$$ where $\gamma_x$ and $\gamma_y$ are the components of $\gamma$. A pseudofold is a subsurface given by such a chart. Note the parameter $d$ depends on the pseudofold angle.

\end{defi}

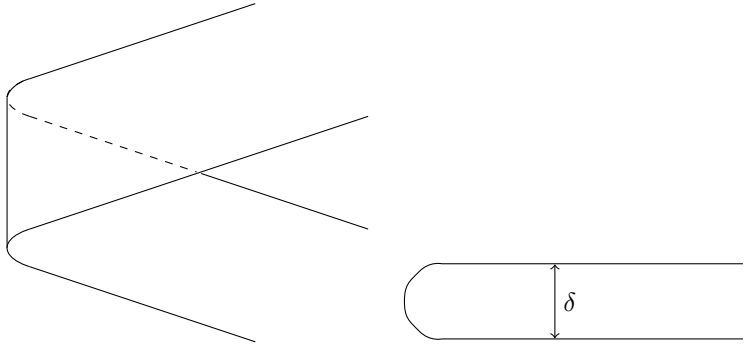
\begin{figure}[!ht]
    \vspace{-20pt}
    \centering
    \begin{tikzpicture}
        \draw (2, -1)--(-1, 0);
        \draw (-1, 0.5)--(3.5, 2);
        \draw (-1, 2.5)--(2, 3.5);
        \draw (1.28, 1.24)--(3.5, 0.5);
        \draw[dashed] (-1, 2)--(1.22, 1.26);
    
        \draw plot [smooth, tension=1.2] coordinates {(-1, 0) (-1.3, 0.25) (-1, 0.5)};
        \draw[dashed] plot [smooth, tension=1.2] coordinates {(-1, 2) (-1.3, 2.25) (-1, 2.5)};
        \draw plot [smooth, tension=1.2] coordinates {(-1.3, 2.25) (-1.2, 2.4) (-1, 2.5)};
        \draw (-1.3, 0.25)--(-1.3, 2.25);
    \end{tikzpicture}
    \quad 
    \begin{tikzpicture}
        \draw (4, 0)--(0, 0);
        \draw plot [smooth, tension=1] coordinates {(0, 0) (-0.3, 0.1) (-0.5, 0.5) (-0.3, 0.9) (0, 1)};
        \draw (0, 1)--(4, 1);
        \draw[<->] (1.5, 0.01)--(1.5, 0.99);
        \node at (1.7, 0.5) {$\delta$};
    \end{tikzpicture}
    \caption{The top and side views of a pseudofold. The side view is simply the graph of $\gamma(\delta, t)$ in the plane.}
    \label{fig:pseudofolds}
\end{figure}

Separate each layer of the folded ribbon knot vertically by some small distance $\delta$ and then connect the layers with pseudofolds. Let the \textit{size} of a pseudofold be the height disparity between the layers it connects. The pseudofolds of the accordion all have size $\delta$. The corresponding to the wrapping have sizes $(m+1)\delta$, $(m+2)\delta$, ..., $(m+n)\delta$, assuming there are $m$ accordion folds and $n$ wrapping folds. We can represent this in a prefold diagram. In this diagram, let green parallelograms represent pseduofolds which replace solid lines, and let purple parallelograms represent pseudofolds which replace dashed lines.

\begin{figure}[!h]
    \centering
    \begin{tikzpicture}
        \node[rectangle, inner sep=0pt, minimum width=0.15cm, minimum height=2.5cm, rotate=-45, fill=Plum!20, draw] at (0, 0) {};
        \node[rectangle, inner sep=0pt, minimum width=0.15cm, minimum height=2.5cm, rotate=-45, fill=Green!20, draw] at (0.5, 0) {};
        \node[rectangle, inner sep=0pt, minimum width=0.15cm, minimum height=2.5cm, rotate=-45, fill=Plum!20, draw] at (1.0, 0) {};
        \node[rectangle, inner sep=0pt, minimum width=0.15cm, minimum height=2.5cm, rotate=-45, fill=Green!20, draw] at (1.5, 0) {};
        \node[rectangle, inner sep=0pt, minimum width=0.15cm, minimum height=2.5cm, rotate=-45, fill=Plum!20, draw] at (2, 0) {};
        \node[rectangle, inner sep=0pt, minimum width=0.15cm, minimum height=2.5cm, rotate=-45, fill=Green!20, draw] at (2.5, 0) {};
        \node[rectangle, inner sep=0pt, minimum width=0.15cm, minimum height=2.5cm, rotate=-45, fill=Plum!20, draw] at (3, 0) {};

        \node[rectangle, inner sep=0pt, minimum width=0.3cm, minimum height=2.6cm, rotate=-45, fill=Green!20, draw] at (-0.6, 0) {};
        \node[rectangle, inner sep=0pt, minimum width=0.3cm, minimum height=2.6cm, rotate=-45, fill=Green!20, draw] at (-0.9, 0) {};
        \node[rectangle, inner sep=0pt, minimum width=0.3cm, minimum height=2.6cm, rotate=-45, fill=Green!20] at (-0.75, 0) {};

        \node[rectangle, inner sep=0pt, minimum width=0.3cm, minimum height=2.6cm, rotate=-45, fill=Green!20, draw] at (-1.6, 0) {};
        \node[rectangle, inner sep=0pt, minimum width=0.3cm, minimum height=2.6cm, rotate=-45, fill=Green!20, draw] at (-1.95, 0) {};
        \node[rectangle, inner sep=0pt, minimum width=0.3cm, minimum height=2.6cm, rotate=-45, fill=Green!20] at (-1.78, 0) {};

        \node[rectangle, inner sep=0pt, minimum width=0.3cm, minimum height=2.6cm, rotate=-45, fill=Green!20, draw] at (-2.65, 0) {};
        \node[rectangle, inner sep=0pt, minimum width=0.3cm, minimum height=2.6cm, rotate=-45, fill=Green!20, draw] at (-3.05, 0) {};
        \node[rectangle, inner sep=0pt, minimum width=0.3cm, minimum height=2.6cm, rotate=-45, fill=Green!20] at (-2.8, 0) {};

        \node[rectangle, inner sep=0pt, minimum width=0.3cm, minimum height=2.6cm, rotate=-45, fill=Green!20, draw] at (-3.75, 0) {};
        \node[rectangle, inner sep=0pt, minimum width=0.3cm, minimum height=2.6cm, rotate=-45, fill=Green!20, draw] at (-4.2, 0) {};
        \node[rectangle, inner sep=0pt, minimum width=0.3cm, minimum height=2.6cm, rotate=-45, fill=Green!20] at (-3.98, 0) {};

        \node[rectangle, minimum height=1.5cm, minimum width=10.5cm, draw] at (-0.5, 0) {};
        \node[rectangle, minimum height=0.3cm, minimum width=10.5cm, fill=white] at (-0.5, 0.905) {};
        \node[rectangle, minimum height=0.3cm, minimum width=10.5cm, fill=white] at (-0.5, -0.905) {};
    \end{tikzpicture}
    \caption{The prefold diagram for the smooth approximation. The base of each parallelogram is $c\sqrt{2}$ times the size of the corresponding pseudofold. $c$ is the same constant used in Definition 5. It depends on the particular curve $\gamma$ used to construct the pseudofolds.}
    \label{fig:enter-label}
\end{figure}
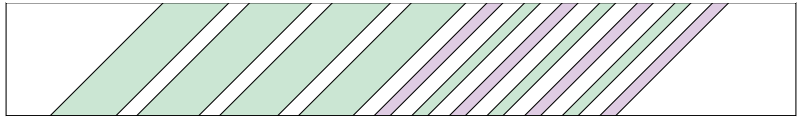

Every error in the approximation is proportional to $\delta$. Thus, as $\delta$ goes to 0, the additional band length and the distance between any particular point on the ribbon and its approximating point on the band go to 0. Furthermore, these bands respect the folding information of $\mathcal{K}_n$, so they are $n$-twist bands.

\end{proof}

\begin{figure}[!hb]
    \vspace{-10pt}
    \centering
    \begin{tikzpicture} 
        \draw (0, 0)--(1, 0);
        \draw (0, 0.1)--(1, 0.1);
        \draw (0, 0.2)--(1, 0.2);
        \draw (0, 0.3)--(1, 0.3);
        \draw (0, 0.4)--(1, 0.4);
        \draw (0, 0.5)--(1, 0.5);
        \draw (0, 0.6)--(1, 0.6);
        \draw (0, 0.7)--(1, 0.7);
        \draw (0, 0.8)--(1, 0.8);
        \draw (0, 0.9)--(1, 0.9);
        \draw (0, 1)--(1, 1);
        \draw (0, 1.1)--(1, 1.1);
        \draw (0, 1.2)--(1, 1.2);
        \draw (0, 1.3)--(1, 1.3);
        \draw (0, 1.4)--(1, 1.4);
        \draw (0, 1.5)--(1, 1.5);
        \draw (0, 1.6)--(1, 1.6);

        \draw (0,0.5) arc [start angle=90, end angle=270, radius=0.05];
        \draw (1,0.5) arc [start angle=-90, end angle=90, radius=0.05];
        \draw (0,0.7) arc [start angle=90, end angle=270, radius=0.05];
        \draw (1,0.7) arc [start angle=-90, end angle=90, radius=0.05];
        \draw (0,0.9) arc [start angle=90, end angle=270, radius=0.05];
        \draw (1,0.9) arc [start angle=-90, end angle=90, radius=0.05];
        \draw (0,1.1) arc [start angle=90, end angle=270, radius=0.05];
        \draw (1,1.1) arc [start angle=-90, end angle=90, radius=0.05];
        
        \draw (0,1.2) arc [start angle=90, end angle=270, radius=0.45];
        \draw (1,0.3) arc [start angle=-90, end angle=90, radius=0.5];
        \draw (0,1.3) arc [start angle=90, end angle=270, radius=0.55];
        \draw (1,0.2) arc [start angle=-90, end angle=90, radius=0.6];
        \draw (0,1.4) arc [start angle=90, end angle=270, radius=0.65];
        \draw (1,0.1) arc [start angle=-90, end angle=90, radius=0.7];
        \draw (0,1.5) arc [start angle=90, end angle=270, radius=0.75];
        \draw (1,0) arc [start angle=-90, end angle=90, radius=0.8];
    \end{tikzpicture}
    \caption{A side view of the complicated part of a smooth approximation}
    \label{fig:enter-label}
\end{figure}
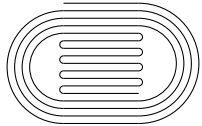

\FloatBarrier

\newpage

\section{Parity}

The construction so far applies to both M\"obius bands and annuli. Which one is constructed comes down to how the ends of the band are connected to one another. Letting one side be colored red and the other blue, a M\"obius band is obtained from taping red to blue, while an annulus is obtained by taping red to red. Below are two reasonably efficient ways to reconnect the ends for each type of band. 

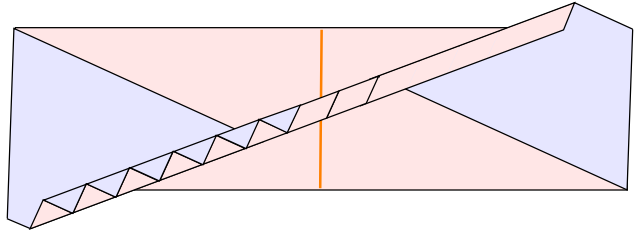
\begin{figure}[!h]
    \vspace{20pt}
    \hspace{110pt}
    \begin{tikzpicture}[transform canvas={rotate=20.5, scale=1.2}]
        \node[rectangle, minimum width=6.78cm, minimum height=1.8cm, rotate=-20.5, fill=red!10, draw] at (1.39, -1.02) {};
    
        \node[isosceles triangle, isosceles triangle apex angle=45, minimum height=2.8cm, rotate=-22.5, fill=blue!10, draw] at (-1.24, -0.35) {};
        \node[isosceles triangle, isosceles triangle apex angle=45, minimum height=2.8cm, rotate=157.5, fill=blue!10, draw] at (4, -1.69) {};
        
        \node[isosceles triangle, isosceles triangle apex angle=90, rotate=225, fill=white] at (-2.4, -1) {};
        \draw (-2.3, -0.94)--(-2.1, -1.13);
        \node[isosceles triangle, isosceles triangle apex angle=90, rotate=45, fill=white, minimum height=0.7cm] at (5, -1.1) {};
        \draw (4.43, -0.9)--(4.93, -1.40);

        \draw[orange, thick] (1.07, -1.84)--(1.7, -0.2);
        
        \node[trapezium, trapezium left angle=45, trapezium right angle=135, fill=red!10] at (2.25, -1.02) {};
        \node[trapezium, trapezium left angle=45, trapezium right angle=135, fill=red!10] at (2.72, -1.02) {};
        \node[trapezium, trapezium left angle=45, trapezium right angle=135, fill=red!10] at (3.19, -1.02) {};
        \node[trapezium, trapezium left angle=45, trapezium right angle=135, fill=red!10] at (3.66, -1.02) {};
        \node[trapezium, trapezium left angle=45, trapezium right angle=135, fill=red!10] at (4.08, -1.02) {};
        \draw (2, -0.90)--(4.43, -0.90)--(4.21, -1.14)--(1.8, -1.14);
        
        \node[trapezium, trapezium left angle=45, trapezium right angle=135, fill=red!10, draw] at (1.76, -1.02) {};
        \node[trapezium, trapezium left angle=45, trapezium right angle=135, fill=red!10, draw] at (1.29, -1.02) {};
    
        \node[isosceles triangle, inner sep=0pt, minimum height=0.25cm, isosceles triangle apex angle=90, rotate=90, fill=red!10, draw] at (0.7, -1.033) {};
        \node[isosceles triangle, inner sep=0pt, minimum height=0.25cm, isosceles triangle apex angle=90, rotate=90, fill=red!10, draw] at (0.195, -1.033) {};
        \node[isosceles triangle, inner sep=0pt, minimum height=0.25cm, isosceles triangle apex angle=90, rotate=90, fill=red!10, draw] at (-0.312, -1.033) {};
        \node[isosceles triangle, inner sep=0pt, minimum height=0.25cm, isosceles triangle apex angle=90, rotate=90, fill=red!10, draw] at (-0.83, -1.033) {};
        \node[isosceles triangle, inner sep=0pt, minimum height=0.25cm, isosceles triangle apex angle=90, rotate=90, fill=red!10, draw] at (-1.34, -1.033) {};
        \node[isosceles triangle, inner sep=0pt, minimum height=0.25cm, isosceles triangle apex angle=90, rotate=90, fill=red!10, draw] at (-1.85, -1.033) {};

        \node[isosceles triangle, inner sep=0pt, minimum height=0.24cm, isosceles triangle apex angle=90, rotate=269, fill=blue!10, draw] at (0.945, -0.995) {};
        \node[isosceles triangle, inner sep=0pt, minimum height=0.25cm, isosceles triangle apex angle=90, rotate=270, fill=blue!10, draw] at (0.455, -0.995) {};
        \node[isosceles triangle, inner sep=0pt, minimum height=0.25cm, isosceles triangle apex angle=90, rotate=270, fill=blue!10, draw] at (-0.055, -0.995) {};
        \node[isosceles triangle, inner sep=0pt, minimum height=0.25cm, isosceles triangle apex angle=90, rotate=270, fill=blue!10, draw] at (-0.57, -0.995) {};
        \node[isosceles triangle, inner sep=0pt, minimum height=0.25cm, isosceles triangle apex angle=90, rotate=270, fill=blue!10, draw] at (-1.085, -0.995) {};
        \node[isosceles triangle, inner sep=0pt, minimum height=0.25cm, isosceles triangle apex angle=90, rotate=270, fill=blue!10, draw] at (-1.595, -0.995) {};
    \end{tikzpicture}
    \vspace{60pt}
    
    \caption{A fully constructed paper annulus. The orange line in the middle indicates where the strip of paper is taped/glued to itself. The fact that there is the same color (red) on each side of the line corresponds to the fact that this is an annulus, not a M\"obius Band.}
    \label{fig:even-connected}
\end{figure}

\begin{figure}[!h]
    \vspace{40pt}
    \hspace{110pt}
    \begin{tikzpicture}[transform canvas={scale=1.2}]
        \node[rectangle, minimum height=2cm, minimum width=3cm, fill=red!10, draw] at (0, 0) {};
        \node[rectangle, trapezium left angle=90, minimum height=2cm, minimum width=3cm, fill=blue!10, draw] at (3, 0) {};
        
        \node[isosceles triangle, isosceles triangle apex angle=45, minimum height=2.8cm, rotate=-22.5, fill=blue!10, draw] at (-1.24, -0.35) {};
        \node[isosceles triangle, isosceles triangle apex angle=45, minimum height=2.8cm, rotate=202.5, fill=red!10, draw] at (4.26, -0.35) {};

        \draw[orange, thick] (1.5, -1)--(1.5, 1);
        
        \node[isosceles triangle, isosceles triangle apex angle=90, rotate=225, fill=white] at (-2.4, -1) {};
        \draw (-2.3, -0.94)--(-2.1, -1.13);
        
        \node[isosceles triangle, isosceles triangle apex angle=90, rotate=-45, fill=white] at (5.4, -1) {};
        \draw (5.3, -0.94)--(5.1, -1.13);

        \node[trapezium, trapezium left angle=45, trapezium right angle=135, fill=red!10, draw] at (2.25, -1.02) {};
        \node[trapezium, trapezium left angle=45, trapezium right angle=135, fill=red!10, draw] at (1.77, -1.02) {};
        \node[trapezium, trapezium left angle=45, trapezium right angle=135, fill=red!10, draw] at (1.29, -1.02) {};

        \node[isosceles triangle, inner sep=0pt, minimum height=0.25cm, isosceles triangle apex angle=90, rotate=90, fill=red!10, draw] at (0.7, -1.033) {};
        \node[isosceles triangle, inner sep=0pt, minimum height=0.25cm, isosceles triangle apex angle=90, rotate=90, fill=red!10, draw] at (0.19, -1.033) {};
        \node[isosceles triangle, inner sep=0pt, minimum height=0.25cm, isosceles triangle apex angle=90, rotate=90, fill=red!10, draw] at (-0.32, -1.033) {};
        \node[isosceles triangle, inner sep=0pt, minimum height=0.25cm, isosceles triangle apex angle=90, rotate=90, fill=red!10, draw] at (-0.83, -1.033) {};
        \node[isosceles triangle, inner sep=0pt, minimum height=0.25cm, isosceles triangle apex angle=90, rotate=90, fill=red!10, draw] at (-1.34, -1.033) {};
        \node[isosceles triangle, inner sep=0pt, minimum height=0.25cm, isosceles triangle apex angle=90, rotate=90, fill=red!10, draw] at (-1.85, -1.033) {};

        \node[isosceles triangle, inner sep=0pt, minimum height=0.24cm, isosceles triangle apex angle=90, rotate=269, fill=blue!10, draw] at (0.945, -0.995) {};
        \node[isosceles triangle, inner sep=0pt, minimum height=0.25cm, isosceles triangle apex angle=90, rotate=270, fill=blue!10, draw] at (0.455, -0.995) {};
        \node[isosceles triangle, inner sep=0pt, minimum height=0.25cm, isosceles triangle apex angle=90, rotate=270, fill=blue!10, draw] at (-0.055, -0.995) {};
        \node[isosceles triangle, inner sep=0pt, minimum height=0.25cm, isosceles triangle apex angle=90, rotate=270, fill=blue!10, draw] at (-0.57, -0.995) {};
        \node[isosceles triangle, inner sep=0pt, minimum height=0.25cm, isosceles triangle apex angle=90, rotate=270, fill=blue!10, draw] at (-1.085, -0.995) {};
        \node[isosceles triangle, inner sep=0pt, minimum height=0.25cm, isosceles triangle apex angle=90, rotate=270, fill=blue!10, draw] at (-1.595, -0.995) {};
    \end{tikzpicture}
    \vspace{40pt}
    \caption{A fully connected many-twist paper M\"obius band. The orange gluing line has opposite colors on either side of it, indicating that the band has an odd number of half-twists.}
    \label{fig:odd-connection}
\end{figure}
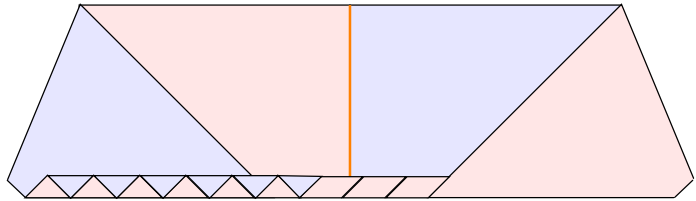

\FloatBarrier

The above constructions give many-twist M\"obius bands with aspect ratio 6.25 and many-twist annuli with aspect ratio 7.45. Thus, we can take $\lambda = 8$ in the Main Theorem. The disparity between cases comes from the fact that the reconnection in the annulus case is less efficient.\footnote{The need for two distinct reconnection methods, and the particular lengths for each type, were worked out together with Luke Briody.}

Note that each reconnection method introduces a handful of folds, and the annulus case includes many center-line crossings. The contributions from the crossings cancel out in pairs, so each reconnection method only contributes a constant amount to the ribbon knot's linking number.

\section{Acknowledgements}

I would like to thank Richard Schwartz and Luke Briody for many helpful discussions around this topic. I am also thankful for the extensive feedback from Schwartz and Elizabeth Denne during the writing process.

\bibliographystyle{plain}
\bibliography{sources}





\end{document}